\documentclass[fullpage,12pt]{amsart}

\makeatletter

\renewcommand{\p@subsection}{}

\renewcommand{\p@subsubsection}{}

\makeatother

\newcommand{\Tr}{{\rm Tr}}
\newcommand{\Fix}{{\rm Fix}}
\newcommand{\Map}{{\rm Map}}

\newcommand{\N}{{\bf N}}
\newcommand{\R}{{\bf R}}
\newcommand{\RP}{{\bf RP}}
\newcommand{\Z}{{\bf Z}}

\title{Euler measure as generalized cardinality}
\author{James Propp, University of Wisconsin}
\address{Department of Mathematics, University of Wisconsin}
\email{propp@@math.harvard.edu}
\date{March 28, 2002}
\begin{document}

\begin{abstract}
Schanuel has pointed out that there are mathematically interesting
categories whose relationship to the ring of integers is analogous to 
the relationship between the category of finite sets and the semi-ring 
of non-negative integers. Such categories are inherently geometrical or
topological, in that the mapping to the ring of integers is a variant 
of Euler characteristic. In these notes, I sketch some ideas that might
be used in further development of a theory along lines suggested by 
Schanuel.
\end{abstract}

\maketitle

\begin{center}
\noindent
\LARGE Foreword 
\normalsize
\bigskip
\end{center}

In this informal article I have gathered together
three memos I wrote in the mid-90s, 
based on conversations with
Scott Axelrod, John Baez, Beifang Chen, 
Timothy Chow, Ezra Getzler, Greg Kuperberg, 
Michael Larsen, Ayelet Lindenstrauss, 
Haynes Miller, Lauren Rose, and Gian-Carlo Rota,
and intended as prologues to further work. 
In the intervening five or six years
my interests have taken me elsewhere,
and I do not expect to return to these topics
anytime soon.
At the same time, I cannot help thinking
that other people might be able to
push these ideas further,
and/or discover that they are more important
(that is, that they are more relevant to
other mathematics) than they currently appear to be.

The first section of this article,
``A proposal for generalizing the Euler characteristic,''
was written in April of 1995.
(Where I've written ``Euler characteristic'' in this section,
the reader should pretend I've written ``Euler measure''.
There are two different ways to generalize Euler's $V-E+F$,
one of which has the nice property of being a homotopy-invariant
and the other of which has the property of being a valuation,
and it seems reasonable to distinguish between them by using
this terminology.)
The second section,
``Negative and fractional cardinalities via generalized polytopes,''
was written in December of 1995,
as an extended abstract for a conference on
Formal Power Series and Algebraic Combinatorics.
The final section,
``Polyhedral sets and combinatorics''
(which should be parsed as ``polyhedral sets-and-combinatorics'',
i.e., polyhedral sets and polyhedral combinatorics),
was written in October of 1996,
to accompany a talk of the same title
that I presented at the Mathematical Sciences Research Institute.
This article concludes with a chapter containing additional comments,
which I have refrained from inserting into the material
that I wrote in the mid-90s.

Another unpublished memo in this vein,
which I wrote at roughly the same time,
has been newly revised by me and will be published
in a special issue of {\it Algebra Universalis\/}
being edited by Joseph Kung
to honor the memory of Gian-Carlo Rota.
It is available over the Web
as an accompaniment to this article.

Rota always encouraged me to pursue my work on Euler measure,
but I was never able to make the sorts of connections
between this work and the broader world of mathematics
that would justify the undertaking.
It's one thing to aspire to do foundational work,
and another thing to have deep insights!
I have often whimsically hoped that someone would create
a journal called ``Definitiones Mathematicae''
that would serve as a haven for interesting definitions
in search of serious theorems that would retroactively justify them.
Lacking such an outlet for my musings,
I have settled for self-publishing these memos
(first on my home-page, and now in the arXiv).
I have made no attempt to remove redundancies between
the three sections of the article.
Also, I have not always included the sort of
bibliographic information that a good scholar should provide
if only for politeness' sake.
If readers of this article have questions,
I'll be happy to try to answer them
(and perhaps include the answers in
later versions of the article).

\begin{center}
\bigskip
\noindent
\LARGE 
I.   A proposal for generalizing the Euler characteristic (1995)
\normalsize
\end{center}

\bigskip

A combinatorialist's fundamental model of a non-negative integer $n$
is a set of $n$ points.  Adding two positive integers corresponds
to taking the disjoint union of two sets; multiplying corresponds
to forming the Cartesian product.

To bring negative numbers into the game, we can follows a suggestion
made by Stephen Schanuel, and replace cardinality by Euler characteristic.
(Note that for finite point-sets, the two notions coincide.)  Thus,
one combinatorial model for the number $-1$ would be a single open
interval (0 vertices and 1 edge yields Euler characteristic $0-1=-1$),
and a model for the negative integer $-n$ would be a disjoint union of 
$n$ open intervals.  Note that our notion of Euler characteristic
is purely combinatorial, and that the sets in question are in general
non-convex and often non-connected; I will call them {\it objects} so 
as not to conflict with established geometric terminology.

In what respect does an open interval $I$ have the properties we expect 
of $-1$?  In the first place, we have relations like $-1+1=0$; but this
is boring.  More interestingly, we have relations like $(-1) \times 
(-1) = 1$: the Cartesian square of an open interval is an open square,
with Euler characteristic $1$.  Even more interestingly, if we
define $X \choose k$ (where $X$ is a topological space and $k$ is a
positive integer) to be the quotient of $X^k$ by the action of the
symmetric group $S_k$, with the part on which $S_k$ does not act freely
removed, then the Euler characteristic of $X \choose k$ is $n \choose k$,
where $n$ is the Euler characteristic of $X$.  For example, if $X$ is
the interval $I$, then $X \choose 2$ is the square $I \times I$ with 
the diagonal removed and with the two resulting pieces identified by
reflection; this new space has Euler characteristic $1$, which is
indeed $-1 \choose 2$.

Most intriguing, however, is the prospect of exponentiation.  In
the case where $Y$ is a finite set of points, we can define $X^Y$
as the set of all functions from $Y$ to $X$, and it will indeed be
the case that $\chi(X^Y) = \chi(X)^{\chi(Y)}$, where $\chi(\cdot)$
denotes the Euler characteristic.  When $Y$ isn't a finite set of
points, but an interval or something even more complicated, then
we clearly won't want to define $X^Y$ as the set of {\it all} maps
from $Y$ to $X$.  But we might define $X^Y$ as the set of all ``nice''
maps from $Y$ to $X$, where niceness is some property or other that
possesses the meta-property that every nice map can be specified by 
a {\it finite number of real-valued parameters} (possibly along with 
some additional combinatorial data).  For instance, the set of nice 
maps from $[0,1]$ to itself could be the set of piecewise-constant 
maps from $[0,1]$ to itself, or the set of all piecewise-linear 
continuous maps from $[0,1]$ to $[0,1]$, or the set of all polynomial 
maps from $[0,1]$ to $[0,1]$.  All of these have the finiteness 
meta-property mentioned above.

Given two objects $X$ and $Y$, and any notion of nice maps from $Y$
to $X$, we might hope to define a generalized Euler characteristic
as the limit, as $n$ goes to infinity, of the (standard) Euler characteristic
of a sequence of objects $P_n$, where the successive $P_n$'s correspond
to sets of nice maps from $Y$ to $X$ with increasingly many (but,
for each $n$, only boundedly many) parameters.  Alternatively, one might 
wish to think of the limit-object $P_\infty$ directly as an 
infinite-dimensional complex having vertices, edges, faces, etc.

\vspace{0.2in}

\noindent
{\bf Example 1:} Consider the the set of piecewise-constant maps from
$(0,1)$ to the two-point set $\{a,b\}$.  We can define the object
$P_n$ as the set of all such maps that are discontinuous at $n$ 
or fewer points in their domain, and we can define $P_\infty$ 
as the direct limit of these objects under the natural inclusion 
maps.  We can then inquire into the behavior of $\chi(P_n)$ as $n
\rightarrow \infty$; we can think of this limit as either
$\lim_{n \rightarrow \infty} \chi(P_n)$ or $\chi(P_0)+\sum_{n=1}^\infty 
\chi(P_n-P_{n-1})$.  I will adopt the latter point of view, and think
of there being vertices, edges, faces, etc.

There are just two vertices (the constant functions).  

What about edges?  These correspond to maps $f$ with a single discontinuity.
If we have a single discontinuity, say at a point $x$, there are three 
things we need to know in order to specify the function $f$: its value 
to the left of $x$, its value to the right of $x$, and its value {\it at} 
$x$.  This would seem to give us 8 possibilities, but in fact we should 
take only 6 (the other 2 correspond to the already-counted constant 
functions).  Each of these 6 combinatorial possibilities ($aab,aba,abb,
baa,bab,bba$) yields an edge of $P_\infty$.

Faces correspond to functions with two discontinuities, $x$ and $y$.  
The set of faces fibers over $I \choose 2$, and each fiber is just a 
string of five letters ($a$'s and $b$'s), in which the first three can't 
all be the same and the last three can't all be the same.  There are 18 
such strings.  (Each of the 6 allowed strings of length three extends 
to 3 strings of length five.)  

Similarly, there are 54 3-cells.

And so on.  Thus the ``Euler characteristic'' is $2-6+18-54+\dots\ $.  
This looks like nonsense, but we can apply Abel summation (or Euler's 
trust-your-pen principle) and assert that this geometric series has 
the value $2/(1-(-3)) = 1/2$.  So $\chi(\{a,b\}^{(0,1)}) = 2^{-1} = 
\chi(\{a,b\})^{\chi((0,1))}$.

\vspace{0.2in}

Note, incidentally, that if we had decided to work in a category in
which our allowed maps were the continuous maps from $(0,1)$ to
$\{a,b\}$, or less trivially the left-continuous maps, we would not
get the answer $1/2$.  So the answer we get seems to be sensitive
to what category we're in.  Still, $2^{-1}$ seems like it should be
the right answer, especially since we can try other experiments in
the category of piecewise-constant maps and ascertain that in many
other cases as well, $\chi(X^Y)=\chi(X)^{\chi(Y)}$.  For instance,
I leave it to you to consider the case $X=Y=I$.

(Fractional Euler characteristics are not in and of themselves novel.
For instance, the infinite-dimensional projective plane $\RP^\infty$
is a 2-to-1 quotient of the infinite-dimensional sphere, which is 
contractible, so it would make sense to define $\chi(\RP^\infty)$ to 
be $1/2$.  Indeed, $\RP^\infty$ has 1 vertex, 1 edge, 1 face, etc., 
so that the preceding ``Eulerian'' method yields Euler characteristic 
$1-1+1-1+\dots=1/(1-(-1))=1/2$.  I'm fairly sure that people have pointed 
this out before --- though I'm not sure who.  Fractional Euler 
characteristics also arise in the theory of group cohomology, and I 
think they show up in orbifold theory as well.  But what's novel here 
is the way in which exponentiation is seen to enter the story.)

Let's try out the piecewise-linear category next.  

\vspace{0.2in}

\noindent
{\bf Example 2:} A linear map $f$ from $[0,1]$ to $(0,1)$ is specified 
by $f(0)$ and $f(1)$, which are arbitrary numbers in $(0,1)$.  So the 
set of piecewise-linear continuous maps from $[0,1]$ to $(0,1)$ with 
{\it no} juncture-points is equivalent to a open square: the Euler 
characteristic is $0-0+1=1$.

What about piecewise-linear continuous maps with a single juncture-point?
Each such map is determined by four numbers: the juncture point $x$
and the values $f(0)$, $f(1)$, and $f(x)$.  Think of the set of
such maps as being fibered over $(0,1) \times (0,1) \times (0,1)$,
corresponding to the choices we make for $x$, $f(0)$, and $f(1)$.
Within each fiber, there is a single forbidden value for $f(x)$,
since we don't want the points $(0,f(0))$, $(x,f(x))$, and $(1,f(1))$ 
to be collinear (we've already counted the maps that are actually linear).
So each fiber is equivalent to $(0,1)$ with a single interior point
removed, and thus has Euler characteristic $-2$.  Multiplying this 
through by $\chi((0,1) \times (0,1) \times (0,1)) = -1$, we find that 
the difference-object $P_1 - P_0$ has Euler characteristic $2$.  

(Note that I'm using facts about Euler characteristic of fiber products;
all this will of course have to be justified, once I figure out exactly
what category I'm in!  For now, though, I'm just trying to get a sense
of what the theory could be like.)

What about piecewise-linear continuous maps with two juncture-points?
Each such map $f$ is determined by the set of juncture points and by
the values of $f$ at 0 and 1, together with its values at the juncture
points themselves.  Think of the first four numbers as determining
a base-space, and the last two as determining a point within a fiber
over that base.  The base space is ${(0,1) \choose 2} \times (0,1) 
\times (0,1)$, which has Euler characteristic 1.  Each fiber is
equivalent (by inclusion-exclusion) to $I \times I - I - I + 1$, 
which has Euler characteristic 4.  Hence $P_2 - P_1$ has Euler
characteristic $1 \times 4 = 4$. 

Continuing in this fashion, one finds that the quasi-Euler characteristic
of the limit-object is $1+2+4+8+\dots$, whose Eulerian value is $-1$.
And sure enough, this is what we should have expected: $\chi((0,1)^{[0,1]})
= \chi((0,1))^{\chi[0,1]} = (-1)^{1} = -1$.  (Thanks to Lauren Rose for
suggesting that I try this example.)

\vspace{0.2in}

\noindent
{\bf Example 3:} Finally, let's consider the set of polynomial maps from 
$[0,1]$ to itself.  For any fixed degree $d$, we can view the set of 
polynomials that take $[0,1]$ into itself as a semialgebraic subset of 
the $d+1$-cube.  This subset is closed, being determined by (uncountably
many) conditions of the form $0 \leq a_0 + a_1 t + \dots + a_d t^d \leq 1$.
It's also bounded (though this requires proof).  Finally, it's contractible, 
because we can take such a polynomial function and multiply it by a 
constant $t$ and send $t$ to zero.  So $P_n$, being a contractible compact
set, has Euler characteristic 1 for all $n$, which would lead us to think 
that the limit-object $P_\infty$ has quasi-Euler characteristic 1 as well.  
And this, too, agrees with our prejudices, since $1 = 1^1 = \chi([0,1])
^{\chi([0,1])}$.

\vspace{0.2in}

This has been a very strange article: no definitions and no theorems!
But I hope the heuristic calculations I've presented are suggestive
of an interesting general theory that might exist.

\begin{center}
\bigskip
\noindent
\LARGE 
II.  Negative and fractional cardinalities via generalized polytopes (1995)
\normalsize
\end{center}

\bigskip
\Large
\noindent
1. Introduction.
\normalsize
\bigskip

For a combinatorialist,
the fundamental significance of the expressions
$n \choose k$ and $k^n$
lies in their interpretation as cardinalities of sets,
specifically,
the set of all $k$-element subsets of an $n$-element set
and the set of $n$-tuples of elements of a $k$-element set.
However, these interpretations are only valid
when $n$ and $k$ are non-negative integers.
In this paper I will describe
an extension of this standard interpretation
that makes sense when $n$ is negative.
Even formulas like ``$\frac12$ choose 2 equals $-\frac18$''
can in some sense be interpreted.

In this enlarged theory,
sets and their cardinalities are replaced
by polytopes and their Euler measures.
The study of Euler measure has its root in Euler's work
on what is now called Euler characteristic,
though the formulation of Euler measure
as an additive set-function
and an explication of its properties
are due largely to Rota, Schanuel, and Chen.
When a polytope is simply a finite collection of points,
its Euler measure is simply its cardinality,
and the standard combinatorial interpretation is recovered
as a special case.

For the present purpose,
I will need to extend the notions of polytope and Euler measure
in two unrelated but compatible directions.
The first extension,
and one that will come as no surprise to those
who know the orbifold notion of Euler characteristic,
is to quotients of polytopes under the free action of a finite group.
The more novel extension is to ``infinite-dimensional polytopes''
of a certain sort.
To assign Euler measure to such objects,
divergent sums must be assigned a notional value
via the physicist's trick of regularization;
more specifically,
an infinite Eulerian sum $V-E+F-\dots$ is interpreted
as the value at $t=-1$ of the holomorphic function
whose Taylor expansion in the vicinity of $t=0$
is $V+Et+Ft^2+\dots$.

Section 2 of this article lays groundwork
by reviewing basic properties of the polyhedral category and Euler measure.
Section 3 discusses quotient polytopes and their role
in providing an interpretation of $n \choose k$
with $n,k \in \Z$, $k \geq 0$.
In Section 4, $\sigma$-polytopes
are introduced and used to interpret $k^n$ with $n,k \in \Z$, $k \geq 1$.
This section also makes it clear why the polyhedral category
is a better setting for generalized combinatorics
than the more familiar topological category.
Section 5 makes a speculative survey of directions
in which the polyhedral approach to combinatorial foundations
can and should be extended,
and Section 6 offers a summary and conclusions.

This article is a preliminary version of a longer article
that I plan to complete before the summer of 1996.
Conversations with
Scott Axelrod, John Baez, Beifang Chen, Ezra Getzler, Greg Kuperberg,
Michael Larsen, Ayelet Lindenstrauss, Haynes Miller, Lauren Rose,
and Gian-Carlo Rota
have helped me clarify my ideas.

\bigskip
\Large
\noindent
2. The polyhedral category and Euler measure.
\normalsize
\bigskip

A {\em polyhedron} is any subset of a Euclidean space $\R^n$
that can be defined through conjunction and disjunction
of a finite number of linear equations and inequalities
involving the $n$ coordinates.
Equivalently, the collection of polyhedra in $\R^n$
is the algebra of sets generated by the (open or closed)
half-spaces of $\R^n$ under union, intersection, and complementation.
For the most part I will focus on bounded polyhedra, or {\em polytopes},
though much of the theory carries over to the unbounded case
(with complications).

The {\em sum} of two disjoint polyhedra in $\R^n$ is defined here
as their union (not their Minkowski sum);
more generally,
the sum of two polyhedra $P,Q \subseteq \R^n$
is $P+Q = P \times \{1\} \cup Q \times \{2\} \subseteq \R^{n+1}$.
The {\em product} of two polyhedra $P \subseteq \R^m$, $Q \subseteq \R^n$
is their usual Cartesian product $P \times Q \subseteq \R^{m+n}$.
A function $f: P \rightarrow Q$ is a {\em polyhedral map}
if its graph (a subset of the polyhedron $P \times Q$)
is also a polyhedron.
Two polyhedra are {\em polyhedrally isomorphic}
if there is a bijective polyhedral map from one to the other.

Examples: The open interval $I = (0,1)$ is a polyhedron in $\R$.
A function from $I$ to $I$ is a polyhedral map
iff it has finitely many ``break-points''
(points of discontinuity or non-differentiability)
and it is linear on the interval between two consecutive break-points.
$I+I$ is polyhedrally isomorphic to $(0,1) \cup (2,3)$ (and to
$(0,\frac12) \cup (\frac12,1)$).  $I \times I$ is the open unit square.

Every polytope in $\R^n$ can be written as a union of finitely many 
(relative-open)
0-cells, 1-cells, $\dots$, and $n$-cells.
The {\em $f$-polynomial} of such a decomposition is defined as
$f(t) = \sum_{i=0}^n f_i t^i \in \N[t]$,
where $f_i$ is the number of $i$-cells in the decomposition.
Two polytopes are polyhedrally isomorphic
iff they have the same dimension and Euler measure.
The {\em dimension} of a polytope $P$ is the largest $d$
for which there exists an injective polyhedral map
from the bounded $d$-dimensional cube into $P$.
The {\em Euler measure} of a $d$-dimensional polytope $P$
is the value of any associated $f$-polynomial at $t=-1$
(the value of this alternating sum is independent
of the decomposition chosen).
We denote the dimension and Euler measure of $P$ 
by $d(P)$ and $\chi(P)$, respectively.
Isomorphism classes of polytopes correspond to
ordered pairs $(d,\chi)$ of integers,
with $d \geq -1$ 
(as usual, we think of the empty set as being $-1$-dimensional);
if $d=-1$, $\chi$ must be 0,
and if $d=0$, $\chi$ must be positive,
while for $d>0$, $\chi$ may be any integer.

If $f_P(t)$ and $f_Q(t)$ are $f$-polynomials
for polytopes $P$ and $Q$,
arising from some specific decompositions,
then these decompositions give rise to
decompositions of $P+Q$ and $P \times Q$
with $f$-polynomials $f_P(t)+f_Q(t)$ and $f_P(t) f_Q(t)$, respectively.
It follows from this that $\chi(\cdot)$,
in addition to being invariant under polyhedral isomorphism,
is finitely additive:
$\chi(P+Q) = \chi(P)+\chi(Q)$.
It also follows that $\chi(\cdot)$ is multiplicative:
$\chi(P \times Q) = \chi(P) \chi(Q)$.

If the polytope $P$ is compact,
then $\chi(P)$ coincides with the Euler characteristic of $P$.
However, unlike Euler characteristic,
Euler measure is not a homotopy invariant;
for instance, the intervals $[0,1]$, $[0,1)$, and $(0,1)$
have Euler measure 1, 0, and $-1$, respectively,
even though they are homotopy-equivalent.

If $f: P \rightarrow \R$ is a piecewise-constant polyhedral map,
with $P$ a polytope,
the {\em Euler integral} $\int_P f \ d\chi$
is defined as
$$\sum_r r \cdot \chi(\{x \in P: f(x)=r\}),$$
where $r$ ranges over the finitely many real numbers in the range of $f$.
(More general versions of the integral can be defined,
but will not be needed here.)
If $f$ is integer-valued, then $\int f \ d\chi$ is an integer;
in particular, for $f$ equal to the indicator function
of the polytope $Q \subseteq P$,
$\int_P f \ d\chi = \chi(Q)$.
The main utility of the Euler integral in this article
is that it facilitates calculation of Euler measure
by way of a Fubini theorem:
for $f: P \times Q \rightarrow \R$,
$$\int_P \left( \int_Q f \ d\chi \right) \: d\chi =
  \int_{P \times Q} f \ d\chi =
  \int_Q \left( \int_P f \ d\chi \right) \: d\chi . $$

\bigskip
\Large
\noindent
3. Binomial coefficients and quotient polytopes.
\normalsize
\bigskip

In laying out some basic notions of quotient polytopes,
I will sidestep the more interesting case of non-free actions,
since they are not germane to my purpose.
Furthermore, I will focus on actions of $G = S_k$,
the symmetric group on $k$ letters,
though everything I say applies to more general free actions
of finite groups.
If a polytope $Q$ is acted on freely by the group $G$,
we can define the associated quotient polytope $G/Q$
as the set of orbits of $Q$ under the action of $G$.
For present intents, it suffices to take $Q$
equal to the set of $k$-tuples
consisting of $k$ distinct points
belonging to some fixed polytope $P$,
with $S$ permuting the $k$ entries;
in this case,
we let $P \choose k$ denote the quotient polytope $Q/S_k$.

One can always find a polytope in $Q$
containing exactly one point in each orbit.
For instance,
define {\em lexicographic ordering} on $P \subseteq \R^n$
in the usual way (with reference to its $n$ coordinates);
then we can represent each orbit in $Q$
by the unique point in that orbit
whose $k$ components are arranged in lexicographically ascending order.
Let us call this the {\it lexicographic representation} of $Q/S_k$.
(E.g., when $P$ is the open interval $I$,
we can represent $I \choose k$ by the set
$\{(x_1,x_2,\dots,x_k) \in I^k: x_1 < x_2 < \dots < x_k\}$,
which, viewed as a subset of $\R^k$,
is just a $k$-dimensional simplex.)
Each polytope that represents $Q/S_k$ 
must have Euler measure $\chi(Q)/k!$.
$\chi(Q)$ itself can be shown to equal
$$\chi(P) \cdot (\chi(P)-1) \cdot (\chi(P)-2) \cdot \cdots (\chi(P)-k+1)$$
by repeated application of the Fubini theorem.
Hence every polytope representing the quotient polytope $P \choose k$
has as its Euler measure the integer $\chi(P) \choose k$
(as given by the standard algebraic definition).

If $P$ is a convex $m$-cell,
then the lexicographic representation of $P \choose k$
can be decomposed in a predictable way
into cells of various dimensions;
for instance, in the case $k=2$, we get
an $m+1$-cell, an $m+2$-cell, \dots, and a $2m$-cell.
A cellular decomposition of the polytope $P$
into convex cells,
with $f$-polynomial $f(t)$,
gives rise to a cellular decomposition of $P \choose 2$
with $f$-polynomial
$$\frac12 \left( [f(t)]^2 - f(t^2) \right) + 
\frac{t}{1-t} \left( f(t) - f(t^2) \right) .$$
Substituting $t=-1$ yields $f(-1) \choose 2$.
A similar situation prevails for $P \choose k$ with $k>2$,
except that the formulas are more complicated.

\bigskip
\Large
\noindent
4. Exponentiation and $\sigma$-polytopes.
\normalsize
\bigskip

Given polytopes $P$ and $Q$,
define $P^Q$ as the set of polyhedral maps $Q \rightarrow P$.
There are natural identifications that can be made
purely at the functional level:
\begin{eqnarray*}
P^{Q+R} & \equiv & P^Q \times P^R \\
P^{Q \times R} & \equiv & (P^Q)^R \\
(P \times Q)^R & \equiv & P^R \times Q^R
\end{eqnarray*}
If $Q$ is a finite set,
then every function from $Q$ to $P$ is polyhedral,
so that the polyhedral definition of exponentiation
coincides with the set-theoretic definition in this case;
$P^Q$ is effectively the $|Q|$th Cartesian power of $P$.

However, when $P$ is finite and $Q$ is infinite,
things get more complicated.
Consider, for instance, the case $|P|=2$, $Q=(0,1)$.
There are only two continuous maps from $Q$ to $P$;
this accounts for the unsuitability of the category of topological spaces
and continuous maps
for the purpose of developing a ``generalized combinatoric''
that features exponentiation.
However, there are infinitely many polyhedral maps.
For any sequence
$0 = x_0 < x_1 < \dots < x_k < x_{k+1} = 1$
we obtain a polyhedral map $f:(0,1) \rightarrow P$
by choosing values $a_0, b_1, a_1, b_2, a_2, \dots, b_k, a_k$ in $P$
and defining $f(x) = a_i$ for $x_i < x < x_{i+1}$ and $f(x_i)=b_i$;
indeed, this representation is unique
if one stipulates that the $x_i$'s are genuine points
of discontinuity for $f$,
i.e., for every $1 \leq i \leq k$,
the values $a_{i-1}$, $b_i$, $a_i$
are not all equal to one another.

If $P$ is finite (say $|P|=m$) and $Q=(0,1)$,
the set of polyhedral maps $f:Q \rightarrow P$
with exactly $k$ discontinuities
has a natural realization
in $\R^k \times P^{2k+1}$ as $S \times F$,
where the simplex 
$S = \{(x_1,x_2,\dots,x_k): 0 < x_1 < x_2 < \dots < x_k\} \subset \R^k$
parametrizes the locations of the discontinuities of $f$
and $F \subset P^{2k+1}$ is the set of
sequences $(a_0, b_1, \dots, b_k, a_k)$ 
for which one never has both $b_i$ and $a_i$ equal to $a_{i-1}$.
It is easily seen that $|F|=m(m^2-1)^k$.
Thus the set of all polyhedral maps $Q \rightarrow P$
can be realized as a union of $m$ 0-cells,
$m(m^2-1)$ 1-cells, $m(m^2-1)^2$ 2-cells, etc.

We cannot evaluate the divergent alternating series 
$m-m(m^2-1)+m(m^2-1)^2-\dots$,
but we can assign it a value through ``regularization''.
If we define the {\em $f$-series}
of this infinite collection of cells
in the obvious way, we get 
$$\sum_{k=0}^\infty m(m^2-1)^k t^k = \frac{m}{1-(m^2-1)t} .$$
Evaluating this at $t=-1$ yields $m/(1+(m^2-1)) = m^{-1} = \chi(P)^{\chi(Q)}$.

More generally, suppose $P$ is a zero-dimensional polytope
consisting of $m$ points
and $Q$ is a one-dimensional polytope composed of
$f_0$ vertices and $f_1$ 1-cells
(where $f_0$ and $f_1$ are not determined by $Q$
but $f_0-f_1=\chi(Q)$ is).
Every polyhedral map from $Q$ to $P$ is determined
by $f_0$ polyhedral maps of a single point into $P$
and $f_1$ polyhedral maps from a 1-cell into $P$.
Thus, our cell-stratification for 
the set of polyhedral maps from a 1-cell into $P$,
along with the obvious stratification for 
the set of polyhedral maps from a 0-cell into $P$,
yield a cell-stratification for $P^Q$
whose generating series is the series expansion of
the rational function
$$ \left( m \right) ^{f_0} \left( \frac{m}{1-(m^2-1)t} \right) ^ {f_1} ;$$
evaluating this function at $t=-1$ yields $m^{f_0-f_1}=\chi(P)^{\chi(Q)}$.

The preceding calculation is related to an alternative way of seeing
that the set of polyhedral maps from $(0,1)$ to an $m$-point set $P$
``ought'' to be assigned Euler measure $\frac1m$,
without explicit recourse to regularization.
On a functional level,
$P^{(0,1)} \equiv 
P^{(0,\frac12)} \times P^{\{\frac12\}} \times P^{(\frac12,1)}$; 
so a desire for functoriality would lead us to want
$\chi(P^{(0,1)}) =
\chi(P^{(0,\frac12)}) \chi(P^{\{\frac12\}}) \chi(P^{(\frac12,1)})$.
On the other hand, the polyhedral equivalence of
$(0,\frac12)$, $(\frac12,1)$, and $(0,1)$
would lead us to expect
$\chi(P^{(0,\frac12)}) = \chi(P^{(\frac12,1)}) = \chi(P^{(0,1)})$.
Combining, we get
$\chi(P^{(0,1)}) =
\chi(P^{(0,1)})^2 \chi(P^{\{\frac12\}})$,
so that either $\chi(P^{(0,1)})=0$
or else $\chi(P^{(0,1)})=1/\chi(P^{\{\frac12\}}) = m^{-1}$.
Putting this differently: If we let $X$ denote $P^Q$,
then there is a nice way to decompose $X$ into $m$ copies of $X \times X$,
so any functor from generalized polytopes to real numbers
that respects $+$ and $\times$
would have to take $X$ either to 0 or to $m^{-1}$.

Objects like $P^Q$ can be construed as
special cases of $\sigma$-polytopes.
I define a {\em $\sigma$-polytope} as 
a formal disjoint union of finite-dimensional cells,
involving only finitely many $k$-dimensional cells for any particular $k$.
It may seem that these objects have too little structure ---
for instance, unlike CW-complexes they carry no information about 
which cells are parts of the boundary of which higher-dimensional cells ---
but this extra information is superfluous in the polyhedral category,
since polyhedral maps (unlike continuous maps)
need not respect boundary-relationships between cells.
Isomorphism classes of $\sigma$-polytopes
correspond to elements of the semi-ring
$$\N[[t]]/(t \sim 2t+1, t^2 \sim 2t^2 + t, t^3 \sim 2t^3 + t^2, \dots);$$
the elements of this semi-ring are equivalence classes of power series
with non-negative integer coefficients,
where two such series are equivalent if each can be obtained from the other
by means of a finite sequence of moves,
each of which replaces a monomial $t^k$
by a sum $2t^{k}+t^{k-1}$
or vice versa.
The geometric significance of the relations $t^k \sim 2t^k + t^{k-1}$
is simple:
every $k$-cell can be divided into two $k$-cells
along with a $k-1$-cell separating them.
Two terminating series (i.e., polynomials in $t$)
are equivalent iff they have the same degree
and the same value at $t=-1$.
Two non-terminating series are equivalent
iff they differ by a polynomial that vanishes at $t=-1$;
in this case, they have the same regularized value at $t=-1$
(assuming that they have a regularized value at $t=-1$ in the first place,
which is not always the case).

It is important to note that the equivalence classes in $\N[[t]]$
that constitute the elements of our semi-ring
are {\em not} closed in $\N[[t]]$ 
relative to the usual ``$t$-adic'' topology
on formal power series in $t$.
Thus, the series
$$2+6t+18t^2+54t^3+162t^4+\dots$$
is equivalent to the series
$$4+12t+18t^2+54t^3+162t^4+\dots,$$
which is equivalent to the series
$$4+12t+36t^2+108t^3+162t^4+\dots,$$
and so on; but none of these series is equivalent to the limit series
$$4+12t+36t^2+108t^3+384t^4+\dots,$$
which is in fact double the original series
(and has regularized value 1, rather than $\frac12$, at $t=-1$).

It should also be noted that the various power series of the form
$$ \left( m \right) ^{f_0} \left( \frac{m}{1-(m^2-1)t} \right) ^ {f_1} $$
with $f_0-f_1$ fixed (but $f_0$, $f_1$ themselves varying)
are typically {\em inequivalent} to each other in the semi-ring,
even though the different polynomials $f_0 + f_1 t$ are {\em equivalent}. 
Thus when we raise one element of our semi-ring
to the power of another,
we should not expect to get a single element
but rather a set of elements. 
We may nevertheless hope that all the elements that we obtain
are equivalent in the weaker sense that
they have the same regularized value at $t=-1$.

As an indication of the compatibility between the ideas sketched
in Sections 3 and 4,
we note that if $P$ is a $\sigma$-polytope
with a generating series $f(t)$ (relative to one particular decomposition),
then there is a natural way to build a $\sigma$-polytope $P \choose 2$
whose elements are unordered pairs of points in $P$;
this $\sigma$-polytope acquires a cellular decomposition
with generating series
$$\frac12 \left( [f(t)]^2 - f(t^2) \right) + 
\frac{t}{1-t} \left( f(t) - f(t^2) \right) ;$$
as long as $f(t)$ has finite regularized value at $t=+1$,
the above expression has regularized value $f(-1) \choose 2$ 
at $t=-1$.
Thus, for instance, if $P$ has regularized Euler measure $\frac12$,
$P \choose 2$ will have regularized Euler measure 
``$\frac12$ choose 2'', or $-\frac18$.

\bigskip
\Large
\noindent
5. Broadening the scope.
\normalsize
\bigskip

One direction in which I am currently extending these ideas
is providing analogous interpretations for $n \choose k$ or $k^n$
in the case where $k$, as well as $n$, is permitted to be negative.

In the case of $n \choose k$,
a natural approach to take is to define $P \choose Q$
(for $P$, $Q$ polytopes)
as the set of polyhedral maps from $P$ to $Q$,
modulo polyhedral bijections of $Q$ with itself.
This is equivalent to the set of polyhedral subsets of $P$
that are polyhedrally equivalent to $Q$,
i.e., that have the same dimension and Euler measure as $Q$. 
As a variant,
one may consider the set of all polyhedral subsets of $P$
that have Euler measure $k$,
for some fixed $k$ (with no constraint on the dimension of the subset).

In the case of $k^n$, the road to take is even clearer:
one should try to find some natural stratification
of the set of polyhedral maps from $Q$ to $P$,
and then verify that the regularized value of the $f$-series at $t=-1$
is $\chi(P)^{\chi(Q)}$.

Another thing to try is to move both $P$ and $Q$ beyond the domain
of 1-dimensional polytopes.
Here we quickly encounter the problem that,
although polyhedral dissections of a 1-dimensional polytope
can be parametrized polyhedrally (by the locations of the break-points),
polyhedral dissections of a 2-dimensional polytope
cannot be so parametrized.
Indeed, to parametrize all the ways of splitting a 2-dimensional polytope
into three pieces by cutting it along a line
(yielding one piece on each side of the line and one piece on the line itself),
we really need to look in the Grassmannian that parametrizes
lines in 2-space.
No doubt recent theories of Euler measure on Grassmannians
will be helpful in this endeavor.

Finally, it would be interesting to try to develop a notion
of generalized Euler measure
in a setting more central to modern mathematics.
Specifically, we could look at the set of {\em continuous}
polyhedral (i.e., piecewise-linear) maps from one polytope to another,
and use the same method of decomposition and regularization
to assign this set of maps an Euler measure.
Piecewise-linear maps,
which can be used to approximate continuous maps
arbitrarily closely,
so in some sense the set of continuous polyhedral maps
might serve as a computational surrogate for
the set of all continuous maps.
This would give us a way to define an Euler characteristic
for the set of continuous mappings from
one topological space to another.
While there is no inherent virtue in making a mere definition,
it seems plausible that the ``combinatorial Euler characteristic'' 
arising under this approach might coincide with the 
``analytic Euler characteristic''
obtained from other, more sophisticated approaches, 
such as Morse theory.
Loop spaces are just one example of a setting
in which this approach might bear fruit.

\bigskip

\Large
\noindent
6. Conclusion.
\normalsize
\bigskip

There are clearly limits to what one should expect from a theory
that purports to ``combinatorify'' exponentiation.
After all,
$$2^{2^{2^{-1}}}$$
is transcendental, while
$$(-1)^{1/2}$$
is complex (and double-valued to boot); 
worse still, once $i$ gets admitted to one's domain of discussion,
the expression $$i^i$$ arises, taking on countably many values.
So we should not expect our system to have good closure properties
under exponentiation.

On the other hand, it is clear that one can go at least some distance
towards the goal of interpreting exponentiation
in a quasi-combinatorial way.
The main problem with the current state of the theory,
in my opinion,
is that I can neither give a recipe for a canonical decomposition
of a $\sigma$-polytope $P^Q$
nor prove that the regularized Euler characteristic
is independent of decomposition
over a broad class of decompositions.
Nevertheless,
I have observed that different ways of trying to calculate
regularized Euler characteristics of various $\sigma$-polytopes
lead to the same answer ---
sometimes for trivial reasons but oftentimes not.
In trying to explain why these different ``meaningless'' calculations
give rise to the same answer,
we may be able to build the substratum of meaning on which they rest.

\begin{center}
\bigskip
\noindent
\LARGE 
III. Polyhedral sets and combinatorics (1996)
\normalsize
\end{center}

\bigskip

A {\em closed convex polytope} $P \subset \R^n$ is a set that can be written
as the intersection of finitely many closed half-spaces.  Given $x \in P$, the
{\em local dimension} of $P$ at $x$ is the maximal $k \geq 0$ for which
$\R^n$ contains a $k$-dimensional (affine) subspace $W$ such that $x$ is 
in the $W$-interior of $P \cap W$.  A {\em k-face} of $P$ is a connected 
component of $\{x \in P:\ \mbox{the local dimension of $P$ at $x$ is $k$}\}$.
More generally, a {\em pure k-cell} in $\R^n$ ($0 \leq k \leq n$) is the
non-empty intersection of a finite number of open half-spaces within a
$k$-dimensional (affine) subspace of $\R^n$.  Every $k$-face is a pure
$k$-cell.

Euler-Poincar\'e Theorem: If $P \subset \R^n$ is an $n$-dimensional non-empty 
compact convex polytope, $F_0 - F_1 + F_2 - F_3 + \dots + (-1)^n F_n = 1$,
where $F_i =$ the number of $i$-faces of $P$.  We can prove this by defining
a suitable valuation (additive function) on a large class of subsets of $\R^n$.

A {\em polyhedral set} in $\R^n$ is 
(1) a union of finitely many pure cells; or, equivalently,
(2) a subset of $\R^n$ described by a finite Boolean formula
involving linear equations and inequalities.
(Schanuel calls it a polyhedral set; Morelli calls it a hedral set.)

Hadwiger-Lenz lemma: There exists a function $\chi(\cdot)$ ({\em Euler 
measure} or {\em combinatorial Euler characteristic}) mapping polyhedral 
sets to integers, such that: 
$\: (1)$ $\chi(A \cup B) = \chi(A) + \chi(B)$ for $A \cap B = \phi$;
$\: (2)$ $\chi(A)=1$ if $A$ is a non-empty compact convex polytope; 
$\: (3)$ $\chi(A) = (-1)^k$ if $A$ is a bounded pure $k$-cell.
(This approach appears earlier in work of Jim Lawrence, Peter McMullen,
and Alexander Barvinok; see also the exposition by Gr\"unbaum and Shephard.)

The Euler-Poincar\'e Theorem is an immediate consequence of the lemma.

The $\chi(A)$ constructed below is invariant under homeomorphisms,
and in the case where $A$ is a PL-manifold in $\R^n$, $\chi(A)$
coincides with the ordinary Euler characteristic; however,
$\chi((0,1)) = -1 \neq +1 = \chi([0,1])$, so $\chi$ does not
coincide with the standard (homotopy-invariant) Euler characteristic.

Check: 
$$\chi(([0,3] \times [0,3]) \setminus ((1,2) \times (1,2))) = 0.$$
$$\chi((\mbox{interior of triangle $abc$}) \ \bigcup \ \{a,b,c\}) = 4.$$
(Note that the latter set is not locally compact, so ordinary homological
approaches to Euler characteristic do not apply.)

Rota and Schanuel's proof of the Hadwiger-Lenz lemma 
uses {\em Euler integration}:
If $f:\R^n \rightarrow \Z$ has the property that $f^{-1}(k)$
is polyhedral for all $k$ and empty for all but finitely many
$k \neq 0$, put $\int f \: d\chi = \sum_k k \chi(f^{-1}(k))$.
E.g., $\int 1_A \ d\chi = \chi(A)$ if $A$ is polyhedral.
Less trivially, if $f: \R \rightarrow \R$ with
$$
f(x) = \left\{ \begin{array}{rl}
	0 & \mbox{if $x<0$} \\
	2 & \mbox{if $x=0$} \\
	-1 & \mbox{if $0<x<1$} \\
	1 & \mbox{if $x=1$} \\
	0 & \mbox{if $x>1$} \end{array} \right.
$$
then 
\begin{eqnarray*}
\int f \: d\chi & = & (2)\chi(\{0\}) + (-1)\chi((0,1)) + (1)\chi(\{1\}) \\
& = & (2)(1) + (-1)(-1) + (1)(1) = 4.
\end{eqnarray*}

Strategy of proof: Integrate $n$-dimensional cross-sectional Euler measure 
$\chi_n$ with respect to 1-dimensional Euler measure $\chi_1$ to define 
$(n+1)$-dimensional Euler measure $\chi_{n+1}$.  
E.g.: If $a,b,c$ are the points $(0,0)$, $(0,1)$, and $(1,0)$ in $\R^2$, 
and $A = (\mbox{interior of triangle $abc$}) \ \bigcup \ \{a,b,c\}$, 
then we put $\chi_2 (A) = \int f \: d\chi_1$ where
$f(x) = \chi_1(\{y: (x,y) \in A\}) =$ the function just discussed;
hence we get $\chi_2(A) = 4$, as before.

Outline of the Rota/Schanuel proof:
Step 1: For all $a$ in $\R$, define $\chi_1(\{a\})=1$ and
$\chi_1((a,\infty)) = \chi_1((-\infty,a)) = -1$,
and for all $a<b$ in $\R$ define $\chi_1((a,b)) = -1$.
Extend $\chi_1$ by finite additivity to all
polyhedral subsets of $\R$.
Step $n'$ ($n \geq 1$): Define $\int f \: d\chi_n = \sum_k k \chi_n(f^{-1}(k))$.
Step $n+1$ ($n \geq 1$): Define $\chi_{n+1}(A) = \int \chi_n(\pi^{-1}(x)) 
\ d\chi_1$, where $\pi(x_1,x_2,\dots,x_{n+1}) = x_{n+1}$.
Then verify that properties (1),(2),(3) hold by induction for all $n$.

\vspace{0.2in}

{\em Finite combinatorics} is the study of the category of finite sets,
with regard to the cardinality functor $\#:{\bf FinSet} \rightarrow \N$.
(E.g.: $\#(A \cup B) = \#(A)+\#(B)-\#(A \cap B)$, $\#(A \times B) = \#(A)\#(B)$,
$\#(A^B) = \#(A)^{\#(B)}$.)

{\em ``Polyhedral combinatorics''} is the study of the category of polyhedral 
sets, with regard to the Euler functor $\chi:{\bf PolySet} \rightarrow \Z$.
(E.g.: $\chi(A \cup B) = \chi(A)+\chi(B)-\chi(A \cap B)$, 
$\chi(A \times B) = \chi(A)\chi(B)$.)

A {\em polyhedral map} $f: A \rightarrow B \ \: (A \subset \R^m$,
$B \subset \R^n$ polyhedral sets) is a function whose graph
is a polyhedral set in $\R^{m+n}$.  E.g.: the piecewise-linear
discontinuous function $f:\R \rightarrow \R$ introduced earlier,
or the function $g(x)=|x|$, or the function $1_A$ for any
polyhedral set $A$.  

Just as two bijectively equivalent finite sets have the same
cardinality, two polyhedrally equivalent sets have the same
Euler characteristic.

A {\em polyhedral permutation} is an invertible polyhedral map
$f:A \rightarrow A$ all of whose orbits have finite cardinality;
we define its {\em trace} $\Tr_A \: f$ as
$\chi(\Fix_A \: f)$ where $\Fix_A \: f = \{x \in A: f(x)=x\}$, 
and its {\em parity} ``$(-1)^f$'' as 
$(-1)^{\chi(D)} = \prod_{k=1}^\infty (-1)^{\chi(O_{2k})}$, where 
$D=\{(x,y) \in A \times A: x<y \ \mbox{but} \ f(x)>f(y)\}$ 
under any linear ordering of $A$ whose graph in $A \times A$ 
is polyhedral, and $O_m$ is any polyhedral set containing 
one representative from each orbit of size $m$.  (Note:
$(-1)^{f \circ g} = (-1)^f (-1)^g$.)

Example: Let $P =$ the union of the three edges of an equilateral
triangle, $G =$ the dihedral group of the triangle 
$=\{$id,flip,flip,flip,rot,rot$\}$.  $\Tr_P({\rm id}) = -3$,
$\Tr_P({\rm flip}) = 1$, $\Tr_P({\rm rot}) = 0$.

Conjecture: Let $P \subset \R^n$ be a union of $k$-faces,
and let $G$ be a group of isometries of $\R^n$ sending
$P$ to $P$.  Define $\rho:G \rightarrow \Z$ by
$\rho(g) = (-1)^k \Tr_P \: g$.  Then $\rho$ is a linear character of $G$.
(Stronger conjecture: This is true for any finite group of polyhedral 
permutations.)  Can we prove this by ``finding the $G$-module''?
[Note: The weaker form of the conjecture was proved independently 
by Miller Maley, Bruce Sagan, Richard Stanley, 
John Stembridge, and Dylan Thurston;
see the Postscript at the end of these notes
for Stembridge's version.]

Conjecture (duality): Let $G$ be a group of isometries of $\R^n$
sending the $n$-dimensional compact convex polytope $P$ to itself.
Let $P_k$ denote the union of the $k$-faces of $P$ ($0 \leq k \leq n$),
determining the character $\rho_k$ as above, and let $\rho_{-1}$ be 
the trivial character.  Also define characters $\rho_k^\circ$
associated with the polar polytope $P^\circ$ of $P$. 
Then there exists a dimension-preserving 
involution on the set of characters of irreducible representations 
of $G$ whose extension by linearity to an involution on the set of 
all characters of $G$ exchanges $\rho_k$ and $\rho_{n-1-k}^\circ$ 
for all $-1 \leq k \leq n$.

\vspace{0.2in}

For $P$ a polyhedral set and $k \geq 0$, consider the action of
the symmetric group $S_k$ on $P^k$.  A {\em free orbit} is one
of cardinality $k!$.  Define ${P \choose k}$ as any polyhedral
set in $P^k$ that contains exactly one point in each free orbit
and no other points.  (Such a polyhedral set exists, and all such
sets are polyhedrally isomorphic.)  Cf.\ Morelli's $\lambda$-ring
structure on the set of polytopes.

Example: $P=(0,1) \subset \R$, ${P \choose 2} = \{(x,y): 0<x<y<1\} =$
a pure 2-cell, $\chi({P \choose 2}) = +1 = {-1 \choose 2} =
{\chi(P) \choose 2}$.  Example: $Q=(0,1) \cup (2,3)$,
\begin{eqnarray*}
{Q \choose 3} & = & \ \ \, \{(x,y,z): 0<x<y<z<1\} \\[-1.5ex]
& & \cup \{(x,y,z): 0<x<y<1,\ 2<z<3\} \\
& & \cup \{(x,y,z): 0<x<1,\ 2<y<z<3\} \\
& & \cup \{(x,y,z): 2<x<y<z<3\},
\end{eqnarray*}
$\chi({Q \choose 3}) = -4 = {-2 \choose 3} = {\chi(Q) \choose 3}$.

Theorem: For any polyhedral set $P$, 
$\chi({P \choose k}) = {\chi(P) \choose k}$.

Given a graph $G=(V,E)$ and a polyhedral set $P$, a {\em P-coloring} of $G$
is a map $f:V \rightarrow P$ such that $f(x) \neq f(y)$ for all $\{x,y\} \in E$.

Theorem: The Euler characteristic of the set of $P$-colorings of $G$ equals
the chromatic polynomial of $G$ evaluated at $\chi(P)$.
E.g.: If $P=\R$, with $\chi(P)=-1$, $P$-colorings of $G$ can be interpreted
directly as points in the complement of the graphical sub-arrangement of the 
braid arrangement determined by $G$.  Since every component 
of the complement of this central hyperplane arrangement has
$\chi = (-1)^{\#(V)}$, this is Zaslavsky's theorem.

\vspace{0.2in}

Fix a polyhedral set $P \subset \R$.  A finite subset $S \subset P$
is {\em fabulous} if for all $t,t' \in S \cup \{+\infty,-\infty\}$,
$\chi((P \setminus S) \cap (t,t'))$ is even.
(Motivation: If $P = \{1,2,\dots,n\}$, a subset of $P$ is fabulous
iff its complement in $P$ can be written as a disjoint union of
pairs $\{k,k+1\}$, and the number of fabulous subsets is the $n+1$st
Fibonacci number.) E.g.: If $P = (0,1) \cup (2,3) \cup (4,5) \cup (6,7)$,
the fabulous subsets of $P$ are $\phi$ and all sets $\{x,y\}$ with 
$2<x<3$, $4<y<5$.

Theorem: If $\chi(P)=n$, the set of fabulous subsets of $P$ has
Euler characteristic $(\varphi^{n+1}-\varphi^{-n-1})/\sqrt{5}$
\ (with $\varphi=(1+\sqrt{5})/2$).
(Generalization to other sequences satisfying linear recurrence
relations?)

\vspace{0.2in}

If $\Phi$ is a collection of finite subsets of $\R^n$, let
$\Phi_k = \{\phi \in \Phi: \#(\phi)=k\}$,
where we identify a $k$-element subset of $\R^n$
with a point in $\R^{kn}$ as before.
We define the ``Euler series'' $\sum_{k=0}^\infty \chi(\Phi_k) t^k$.
If this series converges in a neighborhood of $t=0$ so as to give
unique analytic continuation in a neighborhood of $t=1$, we call
the value at $t=1$ the {\em (regularized) Euler characteristic} of $\Phi$.

Example 1: $\Phi =$ the collection of all finite subsets of $P$,
where $\chi(P)=n$.  Then the Euler series is $1 + nt + {n \choose 2} t^2
+ \dots = (1+t)^n \rightarrow 2^n$ as $t \rightarrow 1$.  (E.g.,
if $P=(0,1)$ with Euler characteristic $-1$, our $\Phi$ has
regularized Euler characteristic $2^{-1} = 1/2$.)

More generally, if $\Phi$ is a ``colored'' collection of finite
subsets of $\R^n$, where each element of $\Phi$ has combinatorial
as well as geometric data, define $\Phi_k$ as the union of the 
$k$-sets in $\Phi$, where $k$-sets of distinct combinatorial type 
are regarded as distinct.

Example 2: $\Phi =$ the collection of all finite subsets of $(0,1)$,
$\Phi' =$ the collection of all 2-element subsets $\{A,B\}$ of $\Phi$.
We can view $\{A,B\}$ as $A \cup B$ equipped with a distinguished
non-empty subset (the symmetric difference of $A$ and $B$) along with
a partition of this set into two subsets ($A \setminus B$ and
$B \setminus A$).  The Euler series is $-t+4t^2-13t^3+40t^4-\dots
= -t/(1+t)(1+3t) \rightarrow -1/8 = {1/2 \choose 2}$.  (This
generalizes to evaluation of the chromatic polynomial of a graph
at any rational number.)  Note that this is not the same approach
as I used in my memo ``Negative and fractional cardinalities via 
generalized polytopes,'' which fails for this case.

Example 3: $\Phi =$ the collection of all {\it polyhedral} subsets
of $(0,1)$.  Every polyhedral $P \subset (0,1)$ determines a finite
set in $(0,1)$, namely its set of {\em break-points} (i.e., points 
of discontinuity of the indicator function of the set); represent
$P$ by its set of break-points, along with combinatorial information
concerning what happens at and between break-points and at the left
and right ends of $(0,1)$, vis-a-vis membership in $P$.  The Euler
series is $2-6t+18t^2-54t^3+\dots = 2/(1+3t) \rightarrow 1/2$.
(More generally, if $P$ is 1-dimensional, the collection of
polyhedral subsets of $P$ has regularized Euler characteristic 
$2^{\chi(P)}$.)

Example 4: $\Phi =$ the collection of all polyhedral subsets of $[0,1)$
of Euler characteristic 0.  The coefficient of $(-t)^n$ in the Euler
series is the central coefficient of $(x+1+x^{-1})^n$, so the Euler
series is $1-t+3t^2-7t^3+19t^4-\dots = \frac{1}{\sqrt{1+2t-3t^2}}$,
which blows up near $t=1$.  (Generalization?)

\vspace{0.2in}

Let $\Map(P,Q) =$ the set of polyhedral maps $P \rightarrow Q$.
When $P$ is finite, $\chi(\Map(P,Q)) = \chi(Q)^{\chi(P)}$.
When $P$ is 1-dimensional, one can still stratify $\Map(P,Q)$
by number-of-break-points, and if moreover $Q$ is finite,
$\chi(\Map(P,Q)) = \chi(Q)^{\chi(P)}$.  But what about
$\chi(\Map(P,Q))$ when $P$ and $Q$ are genuinely 1-dimensional
(e.g., $P=Q=(0,1)$)?  What about $\chi(\Map(P,Q))$ when $P$
is $\geq 2$-dimensional?

What is the right framework for looking at these infinite-dimensional
polyhedral sets?  (Homology theory for non-locally-finite spaces?)

What are the connections between the present theory and the algebraic
enumerative approach to Euler characteristic (counting points on
varieties over finite fields)?

\vspace{0.4in}

\noindent
{\bf Postscript}

\vspace{0.2in}

\noindent

John Stembridge writes:

\vspace{0.2in}

I have a proof of Jim's conjecture about group actions on polyhedral sets.

More specifically, let $G$ be a finite group of isometries of $R^n$ 
that permutes a disjoint set of $k$-cells. (A $k$-cell is by definition 
a $k$-dimensional intersection of open half spaces.) For $g \in G$, define
$f(g) =$ {the Euler-measure of the polyhedral set that is fixed 
pointwise by $g$}.

(Recall that the Euler measure of a $j$-cell is $(-1)^j$.)

CLAIM: $(-1)^k * f$ is the character of a representation of $G$.

BTW: We must insist that the cells are disjoint, or there exist
counterexamples.
 
PROOF. Wlog, we can assume that there is just one orbit of $k$-cells.
Fix a $k$-cell $C$, and let $H$ be the subgroup of $G$ that preserves $C$.
Each $g$ in $H$ permutes the vertices of (the closure of) $C$, so in 
particular $H$ fixes the centroid of $C$'s vertices. Taking 
this centroid as our origin,
let $V_C$ denote the vector space spanned by $C$.  $H$ acts as a group of
isometries of $V_C$. If some $g$ in $H$ has a $j$-dimensional space of fixed
points, then the portion of $C$ that is fixed pointwise by $g$ is a $j$-cell, 
and hence has Euler measure $(-1)^j$. On the other hand, the determinant of
a (real) orthogonal transformation of a $k$-dimensional space has
determinant $(-1)^l$, where $l$ denotes the multiplicity of the eigenvalue 
$-1$.  Since the complex eigenvalues occur in conjugate pairs, we have 
$l+j \equiv k \ \mbox{mod 2}$,
so $\det(g) = (-1)^{(k-j)}$. Using standard rules for inducing representations,
it follows that $(-1)^k * f$ is the character obtained by inducing det
from $H$ to $G$. QED.

\begin{center}
\bigskip
\noindent
\LARGE
Afterword (2002)
\normalsize
\end{center}

\bigskip

One direction that might be interesting to explore is the study of
``polyhedral vector spaces'', as a generalization of the notion of
finite-dimensional vector spaces.  An example of such a vector space
would be the space of polyhedral real-valued functions on the
polyhedral set $A$.  Such spaces have bases, and in every case
I've looked at, there is a natural way to view the set of basis
vectors as a polyhedral set, and what is more, the Euler measure
of the basis turns out to be equal to the Euler measure of $A$.
Is there a general theorem here?

Secondly, as an historical aside, I mention that the surprising
formula $\chi(A^{-1})=0$ that holds when $A$ is a 1-dimensional
polyhedral set satisfying $\chi(A)=0$ and $A^{-1}$ is interpreted
as the set of maps from an open interval into $A$, and which is
proved in the companion article ``Exponentiation and Euler measure,''
is reminiscent of an interesting ``mistake'' made by Brahmagupta of 
Multan in his 6th century treatise {\it Brahmasphutasiddantha\/}.
In that work, Brahmagupta stated rules for manipulating zero in
combination with ordinary numbers: $A+0=A$, $A-0=A$, $Ax0=0$, and 
$A/0=0$.  Of course the last of these is wrong under the usual
understanding of division.  But it is amusing to find a context
in which Brahmagupta's postulate makes sense and is correct.

\begin{center}
\bigskip
\noindent
\LARGE
References
\normalsize
\end{center}

\bigskip

\noindent
Beifang Chen, The Gram-Sommerville and Gauss-Bonnet theorems 
and combinatorial geometric measures for noncompact polyhedra,
{\it Advances in Mathematics} {\bf 91} (1992), 269-291.

\vspace{0.2in}

\noindent
Beifang Chen, On the Euler characteristic of finite unions of convex sets,
{\it Discrete and Computational Geometry} {\bf 10} (1993), 79-93.

\vspace{0.2in}

\noindent
Branko Gr\"unbaum and G.C.\ Shephard, ``A new look at Euler's theorem
for polyhedra,'' {\it American Mathematical Monthly} {\bf 101} (1994),
pp.\ 109-128.  See also the discussion of this article on pp.\ 959-962 
of that same volume and in {\it Mathematics Reviews} 96c:52024.  

\vspace{0.2in}

\noindent
Peter McMullen, ``The polytope algebra,'' {\it Advances in Mathematics}
{\bf 78} (1987), 76--130.

\vspace{0.2in}

\noindent
Robert Morelli, ``A theory of polyhedra,'' {\it Advances in Mathematics}
{\bf 97} (1993), 1--73.

\vspace{0.2in}

\noindent
James Propp, ``Exponentiation and Euler measure,''
{\it Algebra Universalis} (to appear); arXiv: {\tt math.CO/0204009/}.

\vspace{0.2in}

\noindent
Gian-Carlo Rota, Introduction to geometric probability, 1986.

\vspace{0.2in}

\noindent
Stephen Schanuel, ``Negative sets have Euler characteristic and
dimension,'' in {\it Proceedings of Category Theory, 1990\/},
Lecture Notes in Mathematics vol.\ 1488, pp.\ 379--385.

\vspace{0.2in}

\noindent
Stephen Schanuel, ``What is the length of a potato?  An introduction
to geometric measure theory,'' in {\it Categories in Continuum
Physics\/} (1986), Lecture Notes in Mathematics vol.\ 1174, pp.\ 118--126.

\end{document}